\documentclass[3p]{elsarticle}

\usepackage{graphicx}
\usepackage{amssymb}
\usepackage{amsthm}
\usepackage{lineno}

\newtheorem{thm}{Theorem}
\newtheorem{lem}{Lemma}
\newdefinition{rmk}{Remark}
\newproof{pf}{Proof}
\newproof{pot}{Proof of Theorem \ref{thm2}}

\begin{document}
\begin{frontmatter}
\title{Parallel dichotomy algorithm for solving tridiagonal SLAEs}
\author{A.V. Terekhov}
\ead{andrew.terekhov@mail.ru}
\address{Institute of
Computational Mathematics and Mathematical Geophysics,
630090,Novosibirsk,Russia}
\address{Budker Institute of Nuclear Physics, 630090, Novosibirsk,
Russia}


\begin{abstract}
A parallel algorithm for solving a series of matrix equations with
a constant tridiagonal matrix and different right-hand sides is
proposed and studied. The process of solving the problem is
represented in two steps. The first preliminary step is fixing
some rows of the inverse matrix of SLAEs. The second step consists
in calculating solutions for all right-hand sides. For reducing
the communication interactions, based on the formulated and proved
main parallel sweep theorem, we propose an original algorithm for
calculating share components of the solution vector. Theoretical
estimates validating the efficiency of the approach for both the
common- and distributed-memory supercomputers are obtained. Direct
and iterative methods of solving a 2D Poisson equation, which
include procedures of tridiagonal matrix inversion, are realized
using the mpi technology. Results of computational experiments on
a multicomputer demonstrate a high efficiency and scalability of
the parallel sweep algorithm.
\end{abstract}

\begin{keyword}
Parallel algorithm \sep Tridiagonal matrix algorithm (TDMA) \sep
Thomas algorithm  \sep Sweep Method \sep Poisson equation \sep
Alternating Direction Method \sep Fourier Method

\PACS 02.60.Dc \sep 02.60.Cb \sep 02.70.Bf \sep 02.70.Hm
\end{keyword}
\end{frontmatter}

\section{Introduction.}
\sloppy

 The progress in numerical methods of solving "complex
problems" is impossible without applying powerful parallel
computer systems. Thus, it is necessary to investigate numerical
algorithms that allow for efficient parallel implementation.

The problem of solving tridiagonal systems of linear algebraic
equations  \cite{Godunov_Rib,Samarski_Nikolaev,Kon2} is one of the
most frequently solved problems in computational mathematics. The
tridiagonal SLAEs arise in three-point approximation of problems
for ordinary differential equations of second order with constant
and variable coefficients, and also in realization of difference
schemes for equations in partial derivatives
\cite{Samarski_Razn,Ianenko,Morton}. As a rule, tridiagonal SLAEs
are solved using various versions of the direct difference
equation method, that is, a sweep method: monotonic, nonmonotonic,
flux and orthogonal
\cite{Godunov_Rib,Samarski_Nikolaev,Thomas,Godunov_ortho,Ilin_Kuznecov}.

Development and improvement of parallel sweep algorithms is of
great interest, which is confirmed by numerous publications
\cite{Kon2,Konovalov,tridiag-exmp1,tridiag-exmp2,Swarztrauber_par,stone_par,wakatani}
concerned with this difficult problem. Analyzing papers dealing
with this topic, we can conclude that presently available parallel
sweep algorithms are insufficiently efficient and, what is more
important, they are insufficiently scalable. The primary cause is
that efficient, in a theoretical aspect, parallel algorithms
realized on different multiprocessor computer systems become
disadvantageous due to the presence of such operations as
communications and synchronizations.

Solving problems by finite-difference methods frequently requires
to solve not one, but a series of tridiagonal SLAEs with different
right-hand sides, the number of problems in the series can reach
thousands. Thus, the problem of designing an efficient parallel
sweep algorithm for solving series of tridiagonal systems of
equations deserves consideration.

In this paper, we propose a new approach to designing a parallel
sweep algorithm for solving a series of tridiagonal SLAEs with a
constant matrix and different right-hand sides. The process of
solving the problem is subdivided into two steps. The first,
preliminary step consists in fixing some rows of the SLAE inverse
matrix by means of a sequential procedure. Then follows
calculation of solutions for all right-hand sides; doing so, for
increasing the algorithm efficiency using the formulated and
proved main parallel sweep theorem, we proposed an original
algorithm for calculating individual components from the solution
vector.

\section{Statement of the problem.}
The series of systems of algebraic linear equations with a
symmetrical constant tridiagonal matrix means
\begin{equation}
\label{main_eq3}
 A {\bf X_{n}}= { \bf F_{n}} ,\quad n=1,...,N.
\end{equation}

$$
A=\left\|%
\begin{array}{cccccc}
  b_1 & a_1 &  &  &  &  \Large 0\\
  a_1 & b_2 & a_2 &  &  &  \\
   & a_2 & b_3 & a_3 &  &  \\
   &  & \ddots & \ddots & \ddots &  \\
   &  &  & a_{n-2} & b_{n-1} & a_{n-1} \\
   0&  &  &  & a_{n-1} & b_n \\
\end{array}%
\right\|
$$, where $N$  is the number of problems in the series.

Assuming that system  (\ref{main_eq3}) is nondegenerate, we aim
for designing a parallel algorithm for solving the problem and
subsequent realizing on a multicomputer
\cite{Voevodin_MPI,Malysh}.

\textbf{Data decomposition.}The computational and communication
complexity of a parallel algorithm, hence, the execution time
depend drastically on the way of decomposition of problem data.
Let us dwell on the problem of mapping the data of problem
(\ref{main_eq3}) onto a set of processor elements (PEs).

Designing algorithms for distributed-memory supercomputers, it is
necessary to take into account the fact that local data (data in
the local memory of the same PE) are accessed much faster than
data on a distant PE. Thus, even during designing the parallel
sweep algorithm, we perform computing such as to minimize
communication interactions by means of increasing the portion of
local calculations (calculations performed with local data).

We ground the proposed parallel sweep algorithm on the following
specification of data distribution between PEs:
\begin{enumerate}
\item Assuming that the number of PEs is  $p$, divide vectors  ${ \bf F_n}$  and ${ \bf X_n}$ into subvectors ${\bf Q_{n_i},\, U_{n_i}}$ as
follows\footnote{The number of elements of a vector  $\bf V$, is
denoted as size $\it{size}\{\bf V\}$.}:

\begin{equation}
{\bf F_n}=\left({ \bf Q_{n_1}}, { \bf Q_{n_2}},...,{\bf Q_{n_p}}
\right)^\mathrm{T}=\left(f^n_1,f^n_2,...,f^n_{\it{size}\{\bf{
F_n}\}-1},f^n_{\it{size}\{{\bf F_n}\}}\right)^\mathrm{T
}\label{decom_f},
\end{equation}

\begin{equation}
{\bf X_n}=\left({\bf U_{n_1}},{\bf U_{n_2}},...,{\bf U_{n_p}}
\right)^\mathrm{T}=\left(x^n_1,x^n_2,...,x^n_{\it{size}\{\bf{X_n}\}-1},x^n_{\it{size}\{{\bf
X_n}\}}\right)^\mathrm{T} \label{decom_x}.
\end{equation}

\item The sizes of  $\bf Q_{n_i}$ and $\bf U_{n_i}$ are chosen under the conditions
$$
\begin{array}{l}
 size\{{\bf Q_{n_i}}\}=size\{{\bf U_{n_i}}\} \geq 2 \quad i=1,...,p\\\\
\sum_{i=1}^{p}\it{size}\{{\bf
Q_{n_i}}\}=\sum_{i=1}^{p}\it{size}\{{\bf
U_{n_i}}\}=\it{size}\{{\bf F_n}\}=\it{size}\{{\bf X_n}\}
\end{array}
$$
\item Demand that the pair of subvectors $\left({\bf Q_{n_i}},{\bf U_{n_i}}\right)$  belong to PE number $i$ .

\item The row of $A$ number $j$  is on the same PE as the pair of elements  $\left(x^n_j,f^n_j\right)$
from (\ref{decom_f}),(\ref{decom_x}).
\end{enumerate}

We should note that the specification of decomposition of  $ \bf
X_n,\:F_{n}$, and $A$ rules out absolutely duplication of the
problem data. Exactly for this distribution we will design the
parallel sweep algorithm for solving of problem (\ref{main_eq3}).

\section{Parallel sweep algorithm}
\subsection{Basic algorithm}

\begin{lem} Let the tridiagonal system of linear equations $A{\bf
X_n}={\bf F_n}$  be divided into subsystems of the form
\begin{equation}
\label{red_system} A_i{\bf U_{n_i}}={\bf Q_{n_i}}, \quad i=1,...,p
,
 \end{equation}

$$
A_i=\left\|%
\begin{array}{cccccc}
  1 &  &  &  &  &  \Large 0\\
  a_{l_i+1} & b_{l_i+2} & a_{l_i+2} &  &  &  \\
   & a_{l_i+2} & b_{l_i+3} & a_{l_i+3} &  &  \\
   &  & \ddots & \ddots & \ddots &  \\
   &  &  & a_{l_i+t_i-2} & b_{l_i+t_i-1} & a_{l_i+t_i-1} \\
   0&  &  &  &  & 1 \\
\end{array}%
\right\|,
$$

$$
t_i=size\{{\bf U_{n_i}}\},\quad  l_i=\sum_{k=1}^{i-1} t_i
$$

according to the proposed approach, and let we know the values of
elements \footnote{The first and last elements of some vector $\bf
V$, are denoted as  $first\{{\bf V}\}$ and $last\{{\bf V}\}$.}
$first\{{\bf U_{n_i}}\}$ $\quad last\{{\bf U_{n_i}}\}$, then
systems (\ref{red_system}) can be solved independently, equality
(\ref{decom_x}) will be fulfilled.
\end{lem}

The proof of the lemma follows from the tridiagonal matrix of the
structure itself.

\textbf{Algorithm 1.} Based on Lemma 1, we can propose the
following parallel sweep algorithm
\begin{enumerate}
\item Decompose initial system  (\ref{main_eq3}) into subsystems of form (\ref{red_system}).
 \item Find solutions in boundary elements   $first\{{\bf U_{n_i}}\} ,\, last\{{\bf U_{n_i}}\}$.
 \item Compute $\bf X_n$ , by solving independently subsystems.
\end{enumerate}

Thus, the parallel sweep algorithm enables to reduce the solution
of problem (\ref{main_eq3}) to $p$ independent subproblems of form
(\ref{red_system}), if values of  $first\{{\bf U_{n_i}}\}$ and
$last\{{\bf U_{n_i}}\}$ are known. However, we have still to solve
the issue of the efficient way of computing "boundary" elements,
i.e., it is necessary to design a parallel algorithm for computing
an individual component of the solution vector $\bf X_n$.

\subsection{Computing arbitrary solution component.}
\begin{lem}

Let $B$ be the symmetrical tridiagonal matrix. Then the value of
the $k$th component of solution vector $\left({\bf Y}\right)_k$ of
the equation $B{\bf Y}={\bf F}$ can be found as

\begin{equation}
\left({\bf Y}\right)_k={ \bf G_k}^\mathrm{T}{\bf F},
\label{eq_orth-l}
\end{equation}
where the vector  $\bf G_k$ is the solution of the following
system of equations
\begin{equation}
 B{ \bf G_k}=\bf e_k \label{b_system},
\end{equation}
$\bf e_k$ is the unit vector.
\end{lem}

\begin{pf} Let $\bf B_{.k}$ be designated by the  $k$th column of $B$,
and $\bf B_{k.}$  by the $k$th row, respectively

$$
B=\left( {\bf B_{.1}}, {\bf B_{.2}},...,{\bf B_{.n}}\right)=\left[
\begin{array}{c} {\bf B_{1.}}\\{\bf B_{2.}}\\...\\{\bf B_{n.}}
\end{array}
\right]
$$
By virtue of the definition of the inverse matrix $BB^{-1}=\bf I$,
the solution of system  (\ref{b_system}) is the  $k$th column of
$B^{-1}$ and the $k$th row (under the condition of matrix
symmetry).

From this follows
\begin{equation}
\begin{array}{ll}
\label{green_v} {\bf G_k}=B^{-1}
\bf{e_k}=\bf{B_{.k}^{-1}}=\left(\bf{B_{k.}^{-1}}\right)^\mathrm{T},\quad
\left(Y\right)_k={\bf B}_{k.}^{-1}{\bf F}={\bf G^\mathrm{T}_kF}.
\end{array}
\end{equation}
Thus, the lemma has been proved.
\end{pf}

Let us represent the algorithm for computing the arbitrary
solution component for the series of tridiagonal equations of form
(\ref{main_eq3}).

\textbf{Algorithm 2}. For calculating  $M$ different components of
the solution vector $\left({\bf X_n}\right)_{k_m},$ $n=1,...,N$;
$m=1,...,M$ from series of equations  (\ref{main_eq3}),
следует:\begin{enumerate}
    \item Find ${\bf G_{k_m}},\, m=1,...,M$, as the solution of  $A{\bf G_{k_m}}=\bf e_{k_m}$.
    \item Define $\left({\bf X_n}\right)_{k_m}$  for the whole series as $$\left({\bf X_n}\right)_{k_m}={\bf G_{k_m}^{\mathrm{T}}}{\bf F_n},\: n=1,...,N, \quad
    m=1,...,M.$$
\end{enumerate}

Thus, Algorithm 2 makes it possible to find a separate component
of the solution vector. It is important that for different $k_m$
the values of $\left({\bf X_n}\right)_{k_m}$ can be calculated
independently. We should note that the vector ${\bf G_{k_m}}$ does
not depend on the right-hand side of (\ref{main_eq3}) , hence, it
may be determined once for all $\bf F_n$.

Let us consider Algorithm 2 in application to the parallel sweep
algorithm. According to Algorithm 1, at the second step it is
required to calculate values of elements $first\{{\bf
U_{n_i}}\},\:$ $ last\{{\bf U_{n_i}}\}$ for all PEs. The specified
data decomposition assumes that only one subvector ${\bf U_{n_i}}$
and one subvector ${\bf Q_{n_i}}$ are placed on a single PE,
therefore, each PE will compute two elements $\left(first\{{\bf
U_{n_i}}\}, last\{{\bf U_{n_i}}\}\right)$ from the solution
vector. Parallel realization of Algorithm 2 entails a great
difficulty, namely, according to (\ref{eq_orth-l}), each PE will
have to perform about $O(size\{{\bf X_{n}}\})$ operations
regardless of the number of involved computational resources.

This causes the problem of modifying Algorithm 2 in such a manner
that the number of operations per PE is about $\displaystyle
O(size\{{\bf U_{n_i}}\})$.

\subsection{The main parallel sweep theorem.}
We will start designing an efficient algorithm for computing
"boundary" elements; for illustration, let us consider the
following boundary-value problem

\begin{equation}
\label{example}
\begin{array}{lll}
\displaystyle \frac{\mathrm{d}^2\varphi}{\mathrm{d}x^2}=-\rho(x);&
\displaystyle \varphi(x_0)=0,\quad \varphi(x_1)=0
\end{array}
\end{equation}

As is known \cite{math_phys}, the solution of problem
(\ref{example}) may be represented in the integral form via the
corresponding Green function \footnote{In this case, we restrict
ourselves only to the fact of its existence.}
\begin{equation}
\label{int_green}
\varphi(x)=\int_{x_0}^{x_1}G(x,s)\rho(s)\mathrm{d}s.
\end{equation}

Let us partition the interval $(x_0,x_1)$ by three points
$\{x_{1/4},x_{1/2},x_{3/4}\}$ and define the right-hand side of
(\ref{example}) as

\begin{equation}
\label{right_hand_side1} \rho(x)=\left\{
\begin{array}{ll}
0, & x_0\leq x \leq x_{1/4}\\ k(x),& x_{1/4}< x < x_{1/2}
\\0,
& x_{1/2}\leq x \leq x_{3/4}\\ 0,& x_{3/4}\leq x \leq x_{1} \\
\end{array}
\right.
\end{equation}

\begin{figure}[htb]
\center \includegraphics[width=0.4\textwidth]{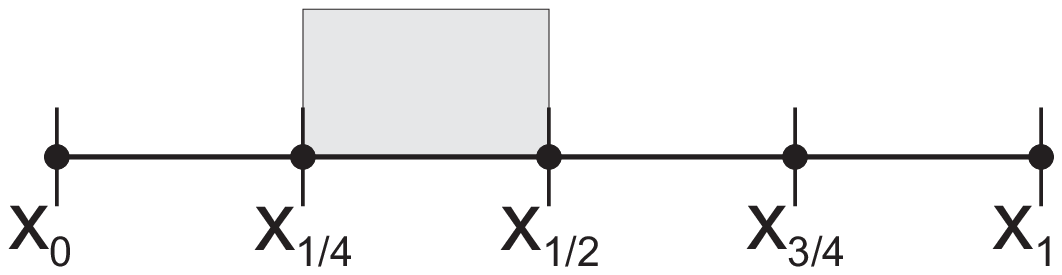}
\caption{}
\end{figure}

According to (\ref{int_green}), the solution at the points of
partitioning may be defined as

\begin{equation}
\begin{array}{ll}
\label{new_method1} \displaystyle
\varphi(x_{1/4})=\int_{x_{1/4}}^{x_{1/2}}G(x_{1/4},s) k
(s)\mathrm{d}s,\quad \displaystyle
\varphi(x_{1/2})=\int_{x_{1/4}}^{x_{1/2}}G(x_{1/2},s) k
(s)\mathrm{d}s,\\
\displaystyle
\varphi(x_{3/4})=\int_{x_{1/4}}^{x_{1/2}}G(x_{3/4},s) k
(s)\mathrm{d}s.
\end{array}
\end{equation}

Another way for finding the solution of equation (\ref{example})
at the point $x_{3/4}$ without calculating the integral of form
(\ref{int_green}), is as follows: since the point  $x_{3/4}$
belongs to the interval  $(x_{1/2},x_{1})$ he solution at it
should satisfy the equation

\begin{equation}
\label{example2} \displaystyle
\frac{\mathrm{d}^2\tilde{\varphi}}{\mathrm{d}x^2}=0,
\end{equation}

with the boundary conditions

\begin{equation}
\label{baundary_example2}
\begin{array}{ll}
\displaystyle
\tilde{\varphi}(x_{1/2})=\int_{x_{1/4}}^{x_{1/2}}G(x_{1/2},s) k
(s)\mathrm{d}s,& \quad \tilde{\varphi}(x_1)=0.
\end{array}
\end{equation}

It is extremely important (from the viewpoint of computation) that
the solution of problem
 (\ref{example2}),(\ref{baundary_example2}) is represented analytically  \cite{sprav_diff}

\begin{equation}
\tilde{\varphi}\left( x \right) =\varphi(x_{1/2}) {\frac { x-{x_1}
{ }}{{ x_{1}}-{x_{1/2}}}}.
\end{equation}

Thus, the solution of  (\ref{example}) with the right-hand side
(\ref{right_hand_side1}) at the points
$\{x_{1/4},x_{1/2},x_{3/4}\}$ is as follows

\begin{equation}
\begin{array}{ll}
\label{new_method2} \displaystyle
\varphi(x_{1/4})=\int_{x_{1/4}}^{x_{1/2}}G(x_{1/4},s) k
(s)\mathrm{d}s,\quad \displaystyle
\varphi(x_{1/2})=\int_{x_{1/4}}^{x_{1/2}}G(x_{1/2},s) k
(s)\mathrm{d}s,\\
\displaystyle \varphi( x_{3/4} ) =\varphi(x_{1/2}) {\frac {
x_{3/4}-{x_1}}{ x_{1}-x_{1/2}}}.
\end{array}
\end{equation}

Comparison of  (\ref{new_method1}) with (\ref{new_method2}) in
their computation complexity shows evident advantage of the latter
because it is required to compute less integrals of form (\ref{int_green}).  \\\\
For the arbitrary function $\rho(x)$ we summarize the obtained
result as a theorem.

\begin{thm} It is required to find the solution of boundary problem
(\ref{example}) at points with the coordinates $\left\{x_i \mid
x_0 <x_i<x_N,\; x_{i} < x_{i+1}, \; i=1,...,\mathrm{N-1}\right\}$.
Then the following identity takes place

\begin{equation}
\label{main_teor1} \varphi(x_i)=\sum_{j=1}^{i}\alpha_j^R {\frac {
x_i- x_N }{x_j-x_N}}+\sum_{j=i+1}^{N} \alpha_j^L {\frac {x_i-
x_0}{x_j-x_0}},\quad i=1,...,N-1,
\end{equation}

\begin{equation}
\begin{array}{cl}
\displaystyle
\alpha_i^R=\int_{x_{i-1}}^{x_i}G(x_i,s)\rho(s)\mathrm{d}s, &
i=1,...,N-1,\\\\ \displaystyle
\alpha_i^L=\int_{x_{i-1}}^{x_i}G(x_{i-1},s)\rho(s)\mathrm{d}s,&
i=2,...,N.
\end{array}
\end{equation}
\end{thm}

Let us formulate and prove the main parallel sweep theorem.
\\\\

\begin{thm} Let we
have the nondegenerate system of linear algebraic equations with
the tridiagonal matrix $A\mathbf{X}=\mathbf{F}$  of dimension $n$.
Then for each solution vector component from the set
\begin{equation}
\label{mnog_xi}
\Omega=\left\{\left( \mathbf{X} \right)_{n_i} \left|1<n_i<n,\;
n_i<n_{i+1},\; i=1,...,p\leq n \right. \right\}
\end{equation} the following identity holds true

\begin{equation}
\label{theor_2}
\left(\mathbf{X}\right)_{n_i}=\sum_{j=1}^{i}\beta_{n_j}^\mathrm{R}
\left(\mathbf{Z}_{n_j}^\mathrm{R}\right)_{n_i}+\sum_{j=i+1}^{p+1}\beta^\mathrm{L}_{n_j}\left(\mathbf{Z}_{n_j}^\mathrm{L}\right)_{n_i},
\end{equation}

\begin{equation}
\label{beta}
\begin{array}{ll}
\displaystyle
\beta_{n_i}^\mathrm{L}=\sum_{j=n_{i-1}}^{n_{i}-1}\mathbf{\textbf(F)}_jA^{-1}_{n_{i-1},j}, & (n_{p+1} = n+1)\\
\displaystyle
\beta_{n_i}^\mathrm{R}=\sum_{j=n_{i-1}}^{n_{i}-1}\mathbf{\textbf(F)}_jA^{-1}_{n_{i},j}, & (n_0 = 1)\\\\
\end{array}
\end{equation}

\begin{equation}
\label{z_vector}
\begin{array}{cc}
B^\mathrm{L}_k\mathbf{Z}^\mathrm{L}_k=\mathbf{e^{\mathrm{L}}},\quad
B^\mathrm{R}_k\mathbf{Z}^\mathrm{R}_k=\mathbf{e^{\mathrm{R}}},
\end{array}
\end{equation}

where

\begin{equation}
\label{matrixbrl}
\begin{array}{cc}
B^\mathrm{L}_k=\left\|%
\begin{array}{cccccc}
  b_1 & c_1 &  &  &  &  \Large 0\\
  a_1 & b_2 & c_2 &  &  &  \\
   & a_2 & b_3 & c_3 &  &  \\
   &  & \ddots & \ddots & \ddots &  \\
   &  &  & a_{k-2} & b_{k-1} & c_{k-1} \\
   0&  &  &  & 0 & 1\\
\end{array}
\right\|, & B^\mathrm{R}_k=\left\|%
\begin{array}{cccccc}
  1 & 0 &  &  &  &  \Large 0\\
  a_k & b_{k+1} & c_{k+1} &  &  &  \\
   & a_{k+1} & b_{k+2} & c_{k+2} &  &  \\
   &  & \ddots & \ddots & \ddots &  \\
   &  &  & a_{n-2} & b_{n-1} & c_{n-1} \\
   0&  &  &  & a_{n-1} & b_{n}\\
\end{array}
\right\|
\end{array}
\end{equation}\

$$
\mathbf{e}^\mathrm{R}=\left(1,0,0,...,0\right)^{\mathrm{T}},\;
\mathbf{e}^\mathrm{L}=\left(0,...,0,0,1\right)^{\mathrm{T}}
$$
\end{thm}

\begin{pf} Let
\begin{equation}
\label{dec}
\begin{array}{cc}
\displaystyle \mathbf{X}=\sum_{j=1}^{p}\mathbf{X}^j,& \quad
\displaystyle \mathbf{F}=\sum_{j=1}^{p}\mathbf{F}^j\\
\end{array}
\end{equation}

$$
\displaystyle A\mathbf{X}^j=\mathbf{F}^j
$$

Define elements of  $\mathbf{F}_j$  as follows

\begin{equation}
\label{define_vectorF}
\left(\mathbf{F}^j\right)_i=\left\{
\begin{array}{ll}
\displaystyle 0, & i < n_{j-1}\\\\
\displaystyle \left(\mathbf{F}\right)_i,& n_{j-1} \leq i < n_j \\\\
\displaystyle 0,& i \geq  n_j
\end{array}
\right.
\end{equation}

Then the solutions  $\left(\mathbf{X}_j\right)_{n_j-1},\;
\left(\mathbf{X}_j\right)_{n_j}$ are unambiguously defined by the
following expressions

\begin{equation}
\begin{array}{c}
\displaystyle \left(\mathbf{X}_j\right)_{n_{j-1}}=A^{-1}_{n_{j-1}\mathbf{.}}\mathbf{F}_j=\sum_{k=n_{j-1}}^{n_j-1}A^{-1}_{n_{j-1},k}\left(\mathbf{F}_j\right)_k=\beta_{n_i}^{\mathrm{L}}\\\\
\displaystyle \displaystyle
\left(\mathbf{X}_j\right)_{n_{j}}=A^{-1}_{n_{j}\mathbf{.}}\mathbf{F}_j=\sum_{k=n_{j-1}}^{n_j-1}A^{-1}_{n_{j},k}\left(\mathbf{F}_j\right)_k=\beta_{n_i}^{\mathrm{R}}
\end{array}
\end{equation}

According to Lemma 1, define   $\left(\mathbf{X}_j\right)_{n_{k}}$
for  $1 \leq k \leq j-2$  as
\begin{equation}
\begin{array}{l}
B^{\mathrm{L}}_{n_j}\mathbf{X}^{\mathrm{L}}_j=\beta_{n_j}^{\mathrm{L}}\mathbf{e}^\mathrm{L},\\
\mathbf{X}_j^\mathrm{L}=\left(x^{\mathrm{L}}_{jn1},...,x^{\mathrm{L}}_{jn_{2}},...,x^{\mathrm{L}}_{jn_{j-2}},...,x^{\mathrm{L}}_{jn_{j-1}}\right)^{\mathrm{T}},
\end{array}
\end{equation}

and for $j+1 \leq k \leq p$ as

\begin{equation}
\begin{array}{l}
 B^{\mathrm{R}}_{n_j}\mathbf{X}^{\mathrm{R}}_j=\beta_{n_j}^{\mathrm{R}}\mathbf{e}^\mathrm{R},\\
\mathbf{X}_j^\mathrm{R}=\left(x^{\mathrm{R}}_{jn_j},...,x^{\mathrm{R}}_{jn_{j+1}},...,x^{\mathrm{R}}_{jn_{p-1}},...,x^{\mathrm{R}}_{jn_{p}}\right)^{\mathrm{T}},
\end{array}
\end{equation}

Denoting
\begin{equation}
\begin{array}{cc}
\mathbf{Z}_{n_i}^\mathrm{R}=\left(B^{\mathrm{R}}_{n_i}\right)^{-1}\mathbf{e}^\mathrm{R},&
\mathbf{Z}_{n_{i-1}}^\mathrm{L}=\left(B^{\mathrm{L}}_{n_{i-1}}\right)^{-1}\mathbf{e}^\mathrm{L},
\end{array}
\end{equation}

we obtain the general formula for computing
$\left(\mathbf{X}_j\right)_{n_i},\; i=1,...,p:$

\begin{equation}
\label{AL} \left(\mathbf{X}_j\right)_{n_i}=\left\{
\begin{array}{l}
\beta_{n_j}^\mathrm{L}\left(\mathbf{Z}_{n_j}^\mathrm{L}\right)_{n_i},\; i > j\\\\
\beta_{n_{j}}^\mathrm{R}\left(\mathbf{Z}_{n_j}^\mathrm{R}\right)_{n_i},\;
i \leq  j
\end{array}\right.
\end{equation}

Substituting (\ref{AL}) into  (\ref{dec}), yields (\ref{theor_2}).
\end{pf}
The theorem has been proved.\\

\begin{rmk} Since the vectors $\mathbf{Z}^\mathrm{R,L}_{n_i}$ do not
depend on the right-hand side of (\ref{main_eq3}), therefore, they
may be defined once for the whole series of problems.
\end{rmk}

\begin{rmk} If  $A=A^\mathrm{T}$, then according to Lemma 2 the
quantities $\beta^\mathrm{R,L}_{n_i}$  may be defined as follows
\end{rmk}

\begin{equation}
\label{beta_symn}
\begin{array}{ll}
\displaystyle
\beta_{n_i}^\mathrm{L}=\sum_{j=n_{i-1}}^{n_{i}-1}\mathbf{\textbf(F)}_j\left(\mathbf{G}_{n_{i-1}}\right)_j,\\
\displaystyle
\beta_{n_i}^\mathrm{R}=\sum_{j=n_{i-1}}^{n_{i}-1}\mathbf{\textbf(F)}_j\left(\mathbf{G}_{n_{i}}\right)_j,
\end{array}
\end{equation}

 \begin{equation}
 \label{green_v}
 A{\bf G_{k_m}}=\bf e_{k_m}.
 \end{equation}
\\\\
Based on Theorem 2, we formulate the algorithm for computing
several components from the solution vector for the series of
tridiagonal equations (\ref{main_eq3}).
\\\\
\textbf{Algorithm 3.}For computing  $M$ different components of
the solution vector $\left({\bf X_n}\right)_{k_m},$ $n=1,...,N$;
$m=1,...,M$ from the series of equations  (\ref{main_eq3}),
follow:\begin{enumerate}
    \item \textbf{Preliminary step.} (Performed once for the whole series of problems).
    \subitem 1.1 Find ${\bf G_{k_m}},\, m=1,...,M$, from the solution of equation  (\ref{green_v}).
    \subitem 1.2 Find ${\bf Z^{\mathrm{R,L}}_{k_m}},\, m=1,...,M$  from
     (\ref{z_vector}).
    \item \textbf{Step of obtaining solutions.} (Performed for each right-hand side F $\mathbf{F}_n,\; n=1,...,N$.)
    \subitem 2.1 Determine $\beta^\mathrm{R,L}_{k_m},\; m=1,...,M$ according to (\ref{beta_symn}) .
    \subitem 2.2 Determine $\left({\bf X_n}\right)_{k_m},\; m=1,...,M$ according to (\ref{theor_2}).

\end{enumerate}

Elementary counting of arithmetic operations at the preliminary
step of Algorithm 3 shows that its realization by formulas
(\ref{z_vector}) and (\ref{green_v}) requires   $\approx 24NP$
operations. For calculating components (\ref{mnog_xi}) by formulas
(\ref{theor_2}),(\ref{beta_symn})  requires  $12N+M^2$ operations.

An important property of Algorithm 3 is that at each of four steps
of the algorithm, calculations for different $k_m$ are
independent. Therefore, the number of arithmetic operations per PE
is $24NM/p$ for the first step and  $(12N+M^2)/p$ for the second,
respectively.

Let us analyze the efficiency of Algorithm \{1,3\} without regard
to communication interactions. As the criterion, we will enter the
speedup

$$
S=T_1/T_p,
$$
where $T_1$ is the number of operations for solving one problem
from series (\ref{main_eq3}) by a sequential sweep algorithm, and
$T_p$ , by Algorithm \{1,3\}. Assuming $p=M$, $T_1=8N$, where $N$
is the number of unknowns and  $T_p=\left(12N+2M^2\right)/p$, we
have \footnote {The obtained estimate (\ref{s_ocenka})  is
conditional and represents rather the qualitative behavior of the
speedup dependence on the number of unknowns and PEs.}

\begin{equation}
\label{s_ocenka}
S=\frac{8Np}{12N+p^2}
\end{equation}

From (\ref{s_ocenka}) it follows that the speedup value increases
monotonically as the number of PEs grows, and then starting from
some  $p>p_0$ decreases monotonically to zero.

Evidently, the minimal time of problem solution is achieved for
the number of PEs

$$p_0=\max_p\left(\frac{8Np}{12N+p^2}\right)=\sqrt{6N},$$ при этом ускорение составит

$$S_{\max}=\frac{\sqrt{6N}}{3}.$$

Thus, the efficiency of parallel Algorithm \{1,3\} for the maximum
possible speedup is  $\approx 30\%$. The remaining  $70\%$
computations fall on "additional" operations for maintaining
parallelism. From (\ref{beta_symn}) and (\ref{theor_2}) it follows
that the volume of these additional computations has order
$O(p^2)$, where $p$ is the number of PEs.

For comparison, that difficulty is also characteristic of the
algorithm proposed in \cite{Konovalov,Paasonen} , where for
computing of $first\{{\bf U_{n_i}}\}$ and $last\{{\bf U_{n_i}}\}$
(in our designation) elements it is necessary to solve the
tridiagonal system of equations with the number of unknowns equal
to the number of PEs. Since the authors propose to calculate the
solution by means of a sequential sweep algorithm version, the
number of additional operations will be of order $O(p)$, but
contrary to (\ref{theor_2}), parallel computing is not allowed.

Thus, in Algorithm \{1,3\} as well as in the algorithm
\cite{Konovalov,Paasonen}, the time of computing $first\{{\bf
U_{n_i}}\}$ and $last\{{\bf U_{n_i}}\}$ elements depends linearly
on the number of PEs.

For increasing the efficiency of Algorithm \{1,3\}, we will task
to reduce the number of arithmetic operations in realizing
formula.

\subsection{Parallel dichotomy algorithm for solving tridiagonal SLAEs.}
It is required to calculate the components of the vector of
solution defined in (\ref{mnog_xi}), it is assumed that
$p=2^{p_0}-1 \leq n,\;p_0>0$.

Let us enter into the consideration the sets
\begin{equation}
\label{omega_i} \displaystyle \Omega_i=\left\{
\left(\mathbf{X}\right)_{n_j} \left| \left( \mathbf{X}
\right)_{n_j} \in \Omega,\; j=2^{ \left \lfloor log_2p \right
\rfloor+1-i}k,\; k=1,...,2^i-1  \right. \right\} \backslash
\left(\bigcup_{j=1}^{i-1} \Omega_{j} \right),
\end{equation} where  $i=1,...,\left \lfloor \log_{2}(p)\right \rfloor+1$.

It is evident that  $$\Omega=\bigcup_{i=1}^{\left \lfloor
\log_{2}(p)\right \rfloor+1}\Omega_i,\quad \quad
\Omega_i\bigcap\Omega_j=\left\{ \emptyset\right\},\; i\neq j$$

\begin{thm} Let the components of the solution vector from the set
$\Omega_j,\; j \geq 1$ and the quanities
$\beta^{\mathrm{R,L}}_{n_i},\;
\mathbf{Z}^{\mathrm{R,L}}_{n_i},\mathbf{G}_{n_j},\;
i=1,...,p,\;j=1,...,p-1$ be determined. Then for all
$\left(\mathbf{X}\right)_{n_i} \in \Omega_{m}, \; j<m \leq \left
\lfloor \log_{2}(p)\right \rfloor+1$, the following identity holds
true

\begin{equation}
\label{theor_3}
\left(\mathbf{X}\right)_{n_i}=\sum_{j=k_1+1}^{i}\beta_{n_j}^\mathrm{R}
\left(\mathbf{Z}_{n_j}^\mathrm{R}\right)_{n_i}+\sum_{j=i+1}^{k_2-1}
\beta^\mathrm{L}_{n_j}\left(\mathbf{Z}_{n_j}^\mathrm{L}\right)_{n_i}
+\delta_{k_1}+\delta_{k_2},
\end{equation}

$$
\begin{array}{l}
 \delta_{k_1}=\left\{
\begin{array}{ll}
0,& k_1=0\\\\
\left(\mathbf{X}\right)_{k_1}\left(\mathbf{Z}_{k_1}^{\mathrm{R}}\right)_{n_i}-\left(\mathbf{G_{k_1+1}}\right)_{k_1}\left(\mathbf{F}\right)_{k_1}\left(\mathbf{Z}_{k_1+1}^{\mathrm{R}}\right)_{n_i},
& k_1>0
\end{array}
  \right.
\\\\
 \delta_{k_2}=\left\{
\begin{array}{ll}
0,& k_2=p+1\\\\
\left(\mathbf{X}\right)_{k_2}\left(\mathbf{Z}_{k_2}^{\mathrm{L}}\right)_{n_i},&
k_2<p+1
\end{array}
  \right.
  \end{array}
$$

where $k_1$ and $\;k_2$ are defined as follows
$$\displaystyle k_1=\min_{t,\;t<k,\;\left(\mathbf{X}\right)_{n_t} \in \left( \Omega_j \bigcup
\{ {X_0}\} \right)}(n_i-n_t),\quad \quad
k_2=\min_{t,\;t>k,\;\left(\mathbf{X}\right)_{n_t}\in
\left(\Omega_j \bigcup \{X_{p+1}\}\right)}(n_t-n_i)$$
\end{thm}
 \begin{pf} Validity of the theorem
follows from the fact that the known components from the set
$\Omega_j$ partition the initial system according to Lemma 1 into
independent subsystems and the solution of each subsystem can be
represented as the sum of general solution of a homogeneous
equation and a partial nonhomogeneous equation
\cite{Gelphond}.\end{pf}

Based on Theorem 3, we formulate the efficient parallel algorithm
for computing separate components from the solution vector
\\\\
\textbf{Algorithm 4. Dichotomy algorithm.} Calculation of  $M$
different components of the solution vector $\left({\bf
X_n}\right)_{k_m},$ $n=1,...,N$, $m=1,...,M$ from series of
equations (\ref{main_eq3}), requires:
   \begin{enumerate}
    \item \textbf{Preliminary step.} (Performed once for the whole series of problems).
    \subitem 1.1 Find ${\bf G_{k_m}},\, m=1,...,M$, from the solution of  (\ref{green_v}).
    \subitem 1.2 Find ${\bf Z^{\mathrm{R,L}}_{k_m}},\, m=1,...,M$  from
     (\ref{z_vector}).
    \item \textbf{Step of obtaining solutions.} (Performed for each right-hand side $\mathbf{F}_n,\; n=1,...,N$.)
    \subitem 2.1 Find $\beta^\mathrm{R,L}_{k_m},\; m=1,...,M$  according to  (\ref{beta_symn}) .
    \subitem 2.2 Calculate in ascending order of index $i=1,...,\left \lfloor \log_{2}(p)\right \rfloor+1$  the components of the solution vector
    $\left({\bf X_n}\right)_{k_m}, \in \Omega_i,$ using (\ref{theor_3}).
\end{enumerate}

\begin{rmk}
Elements belonging to the same set $\Omega_i$, can be calculated
independently.
\end{rmk}

Let us analyze the issue of computational stability of Algorithm
4. We will say that Algorithm 4 is stable if for (\ref{theor_3})
for all $j$
\begin{equation}
\label{stab_cond}
{\left\| \mathbf{Z}^{\mathrm{R,L}}_{n_j} \right\|}_C \leq 1,
\end{equation}
 where

$${\left\| \mathbf{X}\right\|}_C=\max_i\left\{\left|\left(\mathbf{X}\right)_i\right|\right\}$$
\\\\
Let us formulate stability criterion of Algorithm 4.
\begin{thm} Let the matrix $A$ have the diagonal dominance
\cite{Samarski_Nikolaev}
\begin{equation}
\label{diag1}
\left|b_i\right| \geq \left|a_i\right|+\left|c_i\right|,\quad i=2,...,N-1,
\end{equation}
\begin{equation}
\label{diag2} \left|b_1\right| \geq \left|c_1\right|,\quad
\left|b_N\right| \geq \left|a_N\right|,
\end{equation}

and at least in one of inequalities (\ref{diag1}) or (\ref{diag2})
, strict inequality holds, then Algorithm 4 is stable.
\end{thm}
\begin{pf} If the matrix  $A$  has the diagonal dominance, then
obviously the matrices $B^{\mathrm{R,L}}_k$ also have a diagonal
dominance. Following the sweep algorithm \cite{Thomas}, the
solution of system
$B^\mathrm{R}_k\mathbf{Z}^\mathrm{R}_k=\mathbf{e^{\mathrm{R}}}$
may be written as

$$
\begin{array}{ll}
z^{\mathrm{R}}_i=\prod_{i=k}^{i+1}\alpha_i,& i=1,...,k-1,\\\\
\displaystyle \alpha_i=\frac{-c_i}{b_i+a_i\alpha_{i-1}}, & i=2,...,k,\\\\
\displaystyle \alpha_1=-c_1/a_1.
\end{array}
$$

From conditions (\ref{diag1}),(\ref{diag2})follows the inequality
$\left|\alpha_i\right| \leq 1$ \cite{Samarski_Nikolaev}, from
where the following estimate takes place
$$
\left|z^\mathrm{R}_i\right|=\left|\prod_{i=k}^{i+1}\alpha_i\right| \leq 1
$$
We can similarly show that $\left|z^\mathrm{L}_i\right|\leq 1$.
Thus, the presence of diagonal dominance entails stability of
Algorithm 4.
\end{pf}
The theorem has been proved
\\\\
\textit{Remark 5.} Since for calculating all elements from the set
$\Omega_i$ to perform  $O(p)$ arithmetic operations regardless of
the index  $i$, for realizing Algorithm 4, computing $p$
components from the solution vector requires $O(p\log_{2}p)$
operations.
\\\\\
\begin{figure}[htb]
\label{pic_analitical} \center
\includegraphics[width=0.7\textwidth]{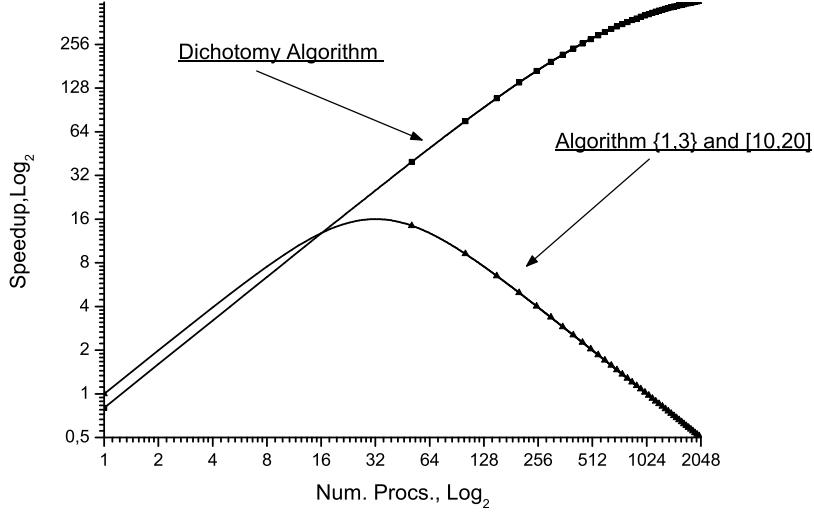}
\caption{Speedup versus the number of PEs.
$size\{\left(\mathbf{X}\right)_n\}=1024$}
\end{figure}

Comparing the dependence of the speedup on the computing time
(Fig. 2) for the dichotomy Algorithm \{1,3\} and algorithm
\cite{Konovalov,Paasonen}, we conclude that the dichotomy
algorithm efficiency for few PEs is comparable with that of
Algorithms \{1,3\} and \cite{Konovalov,Paasonen}. For a great
number of PEs, Algorithms \{1,3\} and \cite{Konovalov,Paasonen}
possess a nearly zero speedup, whereas the dichotomy algorithm
efficiency remains quite high.

\subsection{An example of applying the dichotomy algorithm.}
For illustrating application of the dichotomy algorithm, we will
consider the problem of definition
$$\Omega=\left\{\left(\mathbf{X}\right)_1,
\left(\mathbf{X}\right)_2, \left(\mathbf{X}\right)_3, ,...,
\left(\mathbf{X}\right)_{15}\right\}$$

Let us define the sets $\Omega_i,i=1,...,\left \lfloor log_215
\right \rfloor+1=4$ according to  (\ref{omega_i}).

$$
\begin{array}{l}
\Omega_1=\left\{\left(\mathbf{X}\right)_8\right\},\\ \Omega_2=\left\{\left(\mathbf{X}\right)_4,\left(\mathbf{X}\right)_{12}\right\},\\
\Omega_3=\left\{\left(\mathbf{X}\right)_2,\left(\mathbf{X}\right)_{6},\left(\mathbf{X}\right)_{10},\left(\mathbf{X}\right)_{14}\right\},\\
\Omega_4=\left\{\left(\mathbf{X}\right)_1,\left(\mathbf{X}\right)_{3},\left(\mathbf{X}\right)_{5},\left(\mathbf{X}\right)_{7},\left(\mathbf{X}\right)_{9},\left(\mathbf{X}\right)_{11},\left(\mathbf{X}\right)_{13},\left(\mathbf{X}\right)_{15}\right\}.
\end{array}
$$

Then we calculate at first all elements from $\Omega_1$, and then
$\Omega_2,\Omega_3,\Omega_4$ (Fig. 3).

\begin{figure}[htb]
\label{pic_analitical} \center
\includegraphics[width=0.7\textwidth]{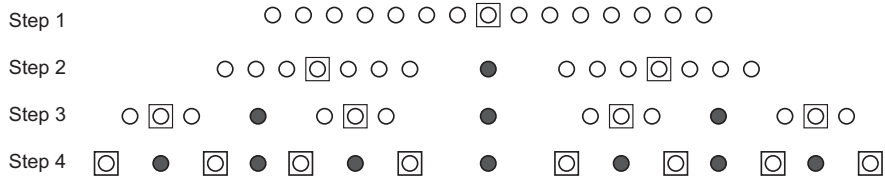} \caption{The
order of computing the elements from the set $\Omega$.}
\end{figure}
Thus, the initial system as a result of calculating at Step 1 the
elements of  $\Omega_1$ is divided into two independent
subproblems, into four independent subproblems at Step 2 after
calculating the elements from  $\Omega_2$ , etc. until calculating
the elements from $\Omega$.

\subsection{Nonsymmetrical matrices}

Until now, it was supposed that the matrix of tridiagonal SLAE
(\ref{main_eq3}) symmetrical. This constraint restricts
considerably the class of problems for which the parallel sweep
algorithm is applicable.

Let us consider problem (\ref{main_eq3}) with the symmetrical
tridiagonal Jacobian matrix whose symmetrical elements have the
same signs ($a_i c_{i-1}\geq0$). In this case, there exists a
similarity transformation with the diagonal matrix $T$ such that
the similar matrix  $\hat{A}=T^{-1} AT$ is symmetrical
\cite{Ilin_fvm}. The elements of the diagonal matrix   $T$ are
defined by the following recurrent relationships

\begin{equation}
\begin{array}{ccc}
T=\rm{diag}\{t_k\},&  t_{k+1}=t_{k}(a_{k+1}/c_{k})^{1/2}, &
k=1,...,N-1
\end{array}
\end{equation}

Thus, the series of SLAEs with the asymmetrical Jacobian matrix
can be solved by the parallel sweep algorithm if the SLAE matrix
is preliminary symmetrized via the similarity transformation.

In the general case, when the tridiagonal matrix is not
symmetrical or cannot be symmetrized, equality (\ref{eq_orth-l})
from Lamma 2 is no longer true. In this case, for determining rows
of the inverse matrix, one can use the explicit representation of
its elements \cite{stone_par,schlegel_inv}

\begin{equation}
\label{x1} A_{ij}^{-1}=\left\{
\begin{array}{rl}
\displaystyle y_iz_j\prod_{k=1}^{j-1}\frac{a_k}{c_{k+1}}, & i \leq j \\
\displaystyle z_iy_j\prod_{k=1}^{j-1}\frac{a_k}{c_{k+1}}, & i \geq
j
\end{array}\right.
\end{equation} ,where\

\begin{equation}
\label{x2} \mathbf{Z}=\left(z_1,z_2,...,z_n\right)^\mathrm{T},
\end{equation}

\begin{equation}
\label{x3}
\begin{array}{cc}
A\mathbf{Y}=\mathbf{e_n},& A\mathbf{Z}=\frac{1}{y_1}\mathbf{e_1}.
\end{array}
\end{equation}

\section{Examples of applying the parallel sweep method.}
For estimating the efficiency of the parallel algorithm of solving
the series of tridiagonal SLAEs we propose, using it as a basis, a
parallel realization of the methods of solving the Poisson
equation. Let us consider the Dirichlet's problem in a rectangle
with the homogeneous boundary conditions $\bar{G}_0=\left\{0\leq
x_{\alpha} \leq l_{\alpha},\,\alpha=1,2 \right\}$

\begin{equation}
\label{poisson}
\begin{array}{c}
\triangle u=-f(x),\quad x=\left(x_1,x_2\right)\in G, \quad \left.
u\right | _\Gamma=0.
\end{array}
\end{equation}

The corresponding difference approximation of second order of
accuracy is

\begin{equation}
\begin{array}{c}
\label{poisson_difference} \Lambda v=-f(x),\quad x\in \omega_h,
\quad \left. v \right | _{\gamma_h}=0\\
\displaystyle (\Lambda y)=\frac{1}{h_1^2}
\left(y_{i+1,j}-2y_{i,j}+
y_{i-1,j}\right)+\frac{1}{h^2_2}\left(y_{i,j+1}-2y_{i,j}+y_{i,j-1}\right),
\end{array}
\end{equation}

where

\begin{equation} \bar \omega_h=\left\{x_i=\left(ih_1,jh_2\right),\quad
i=0,...,N_1,\quad j=0,...,N_2\right\} \label{grid}
\end{equation} is a mesh with steps $h_1$ and  $h_2$, $\gamma_h$ is the mesh boundary.

We will consider the variable separation method (Fourier method)
\cite{Samarski_Nikolaev,Hockney} и Alternating Direction Method
(ADI) \cite{Ianenko,Samarski_Nikolaev,Pissman} with application to
solving problem (\ref{poisson_difference}).

\textbf{a. Variable separation method}. Since the function
$u_{i,j}$ vanishes if обращается в нуль при $j=0$  and $j=N_2$,
and the mesh function  $f_{i,j}$ is given for  $1\leq j\leq
N_2-1$, they may be represented as a series in eigenfunctions of
the difference operator $\Lambda_2 $
\cite{Samarski_andreev,Samarski_Nikolaev}:

\begin{equation}
 \label{operator2}
\left(\Lambda_1y\right)=\frac{y_{i+1,j}-2y_{i,j}+
y_{i-1,j}}{h_1^2}, \quad
\left(\Lambda_2y\right)=\frac{y_{i,j+1}-2y_{i,j}+
y_{i,j-1}}{h_2^2}
\end{equation}

\begin{equation}
\begin{array}{cl}
\displaystyle
u_{i,j}=\sum_{l=1}^{N_2-1}\tilde{u}_i\left(l\right)\sin\left(\frac{\pi
l
j}{N_2}\right) \label{fft_phi}, & \quad 0\leq j \leq N_2,\quad 0\leq i\leq N_1,\\
\displaystyle
f_{i,j}=\sum_{l=1}^{N_2-1}\tilde{f}_i\left(l\right)\sin\left(\frac{\pi
l j}{N_2}\right), & \quad 1\leq j \leq N_2-1, \quad
 1\leq i\leq N_1-1.
\end{array}
\end{equation}

Substituting  (\ref{fft_phi})  into (\ref{poisson_difference})
yields

\begin{equation}
\sum_{l=1}^{N_2-1}\left\{h_1^{-2}\left[\tilde{u}_{i+1}(l)-2\tilde{u}_{i}(l)+\tilde{u}_{i-1}(l)\right]-4h^{-2}_2\tilde{u}_i(l)\sin^2\frac{\pi
l}{2N_2}+\tilde{f}_i(l)\right\}\sin\frac{\pi l j}{N_2}=0
\end{equation}

From this, due to orthogonality of the eigenfunctions
\cite{Samarski_Nikolaev},  the amplitudes of harmonics of the
potential $\tilde{u}_{i}(l),\; l=1,...,N_2-1$ can be defined as
the solution of the following system of equations

\begin{equation}
\begin{array}{lr}
\displaystyle
\tilde{u}_{i+1}(l)-\left(2+4\frac{h_1^2}{h_2^2}\sin^2\frac{\pi
l}{2
N_2}\right)\tilde{u}_i(l)+\tilde{u}_{i-1}(l)=-h_1^2\tilde{f}_{i}(l), & i=1,...,N_1-1,\\\\
\tilde{u}_{0}(l)=\tilde{u}_{N_1}(l)=0. \label{SLAU_FFT}
\end{array}
\end{equation}

The sums (\ref{fft_phi}) should be evidently computed using the
fast discrete Fourier transform \cite{Samarski_Nikolaev,Ilin_fvm},
and for finding the solutions from the series of equations
(\ref{SLAU_FFT}), we should use the sweep method.

Let us dwell on some aspects of realizing the parallel sweep
algorithm within the scope of the variable separation method.

One of the constraints on the parallel sweep algorithm is that all
SLAEs from the series of problems (\ref{main_eq3}) contain the
same fixed matrix. The tridiagonal matrices from (\ref{SLAU_FFT})
have the form

\begin{equation}
\label{series_fft} B_l=(T-d_lI),\quad
d_l=4\frac{h_1^2}{h_2^2}\sin^2\frac{\pi l}{2N_2} , \quad
l=1,...,N_2-1,\quad
\end{equation}

\begin{equation}
{
 \scriptstyle T=\left\|%
\begin{array}{cccc}

  \scriptstyle -2 &  \scriptstyle1 &  &     \scriptstyle\Large 0\\
    \scriptstyle1& \scriptstyle -2 & \scriptstyle 1 &    \\
      & \scriptstyle \ddots &  \scriptstyle\ddots &    \\
       \scriptstyle0&  &    \scriptstyle1 &  \scriptstyle -2 \\
\end{array}%
\right\|}
\end{equation}

It is evident that for (\ref{series_fft}) the condition of matrix
constancy for all right-hand sides is not fulfilled, hence, in
this formulation, problem (\ref{SLAU_FFT})  cannot be solved
efficiently by the proposed algorithm. However, we can extend
Algorithm 1 for solving the series of Poisson equations on the
fixed mesh

\begin{equation}
\triangle u=-f_n(x) ,\quad n=1,...,N.
\end{equation}

In this case, it is required to solve the following problem

\begin{equation}
\label{poisson_end} B_lu_n(l)=g_n(l),\quad l=1,...,N_2-1,\quad
n=1,...,N.
\end{equation}

The set of equations (\ref{poisson_end})  may be considered as a
set of problems of form (\ref{main_eq3}) for the fixed $l$.

\textbf{b. Alternating Direction Method} -- belongs to the class
of methods based on the concept of fixing. The solution of
stationary  problem (\ref{poisson})  is found as the limit  $t
\rightarrow \infty$ of solution of the following unstationary
problem

\begin{equation}
\label{parabolic} \frac{\partial u}{\partial t}=\Delta u-f
\end{equation}  with the same boundary conditions.

Let us consider the Peaceman-Rachford scheme know also an ADI
method \cite{Ianenko,Pissman}. For this purpose, we represent the
2D difference Laplace operator as the sum of two operators
$\Lambda=\Lambda_1+\Lambda_2$ (\ref{operator2}). Then the
iterative process of the ADI method for problem (\ref{parabolic})
has the form

\begin{equation}
\frac{u^{n+1/2}-u^{n}}{\tau^{(1)}_n}=\Lambda_1u^{n+1/2}+\Lambda_2u^n-f,
\label{adi_iter}
\end{equation}

\begin{equation}
\frac{u^{n+1}-u^{n+1/2}}{\tau^{(2)}_n}=\Lambda_1u^{n+1/2}+\Lambda_2u^{n+1}-f.
\label{adi_iter2}
\end{equation}

The iterative parameters  ${\tau^{(1)}_k},\;{\tau^{(2)}_k}$ should
be chosen from the condition of minimum number of iterations. The
problem of choosing the optimal parameters is comprehensively
described, e.g., in
\cite{Samarski_Nikolaev,Samarski_Razn,Wachspress_ADI,Morton}.

Let us consider a case when the region $\bar{G}$ is a square with
the side $l=l_1=l_2$ and the mesh $\bar{\omega}$ is uniform with
$N_1=N_2=N$. Then in order that under any initial approximation
$u_0$ the norm of initial error to be decreased $1/\varepsilon$
times $$
\parallel u_n-u\parallel_D \leq \varepsilon \parallel u_0-u
\parallel _D
$$

the number of iterations $n$ must satisfy the condition

\begin{equation}
\label{ocenka1} n\geq
n_0(\varepsilon)=0.2\ln\left(4N/\pi\right)\ln\left(4/\varepsilon\right).
\end{equation}

Taking into account the fact that the sequence of optimal
parameters $\tau^{(1)}_k,\tau^{(2)}_k,\quad k=1,...,n_0$  ,  is
cyclic and the series of SLAEs (\ref{adi_iter}),(\ref{adi_iter2})
for the fixed $n$  includes the constant matrix

$$
\label{series_adi1}
C^{(1)}_n=\left(T-\frac{h_1^2}{\tau^{(1)}_n}I\right),\quad
C^{(2)}_n=\left(T-\frac{h_2^2}{\tau^{(2)}_n}I\right)
$$

we conclude that at the preliminary step, it is sufficient to
solve merely $n_0$  tridiagonal SLAEs. It should be noted that the
value $n_0$ is much less than the total number of equations whose
solutions have to be found for achieving the desired accuracy.

\section{Computational experiments.}
As we have already mentioned, the parallel algorithms that are
efficient from the theoretical viewpoint, when realized on
supercomputers, may not ensure the expected reduction of the
computation time. The primary reason is that in analyzing the
efficiency of a particular algorithm, it is not easy to take into
account all peculiarities of computer systems (memory operation,
network throughput and latency, etc.). Thus, numerical experiments
with model formulations of problems are an important stage of
investigating parallel algorithms.

As the model problem we considered the Dirichlet problem for the
Poisson equation

\begin{equation}
\label{poisson_test}
\begin{array}{c}
\triangle u=-8\pi^2\sin(2 \pi x)\sin(2 \pi y),\quad
x=\left(x_1,x_2\right)\in G, \quad \left. u\right | _\Gamma=0.
\end{array}
\end{equation}
$$\bar{G}=\left\{0\leq x_{\alpha} \leq 1,\,\alpha=1,2
\right\}$$

For solving problem (\ref{poisson_test})  in the Fortran-90
language using the MPI technology we realized Fourier and ADI
methods. The tridiagonal matrices were inverted by a parallel
dichotomy algorithm. Equation (\ref{poisson_test})  was
approximated on uniform mesh (\ref{grid}) with $N_1=N_2=2^k$
nodes. For the ADI method, the value of prescribed accuracy
$\varepsilon$  was $10^{-5}$.

Figure 3a represents calculation domain decomposition for the
Fourier method. Solution of the tridiagonal systems of equations
was performed in the direction $k_2$, and the Fourier transform
was done in the direction $k_1$. For the ADI method we chose a
decomposition like a lattice (Fig. 3b) because the ADI method
requires solution of tridiagonal SLAES in the directions $x$ and
$y$.

The computing time was estimated as the average time of solving
one problem like (\ref{poisson_test}) from a series of 100
problems

$$
T_{avr}=\frac{\sum_{i=1}^{100}T^i}{100},
$$
and the speedup time was calculated from the formula
$$S_{avr}=\frac{T_{avr}}{T_1},$$ $T^{i}$ is the time of solving problem
(\ref{poisson_test}) by the parallel algorithm, and  $T_1$ -- the
sequential algorithm.

\begin{figure}[!h]
\label{main_pic23}
\begin{center}
\includegraphics[width=0.5\textwidth]{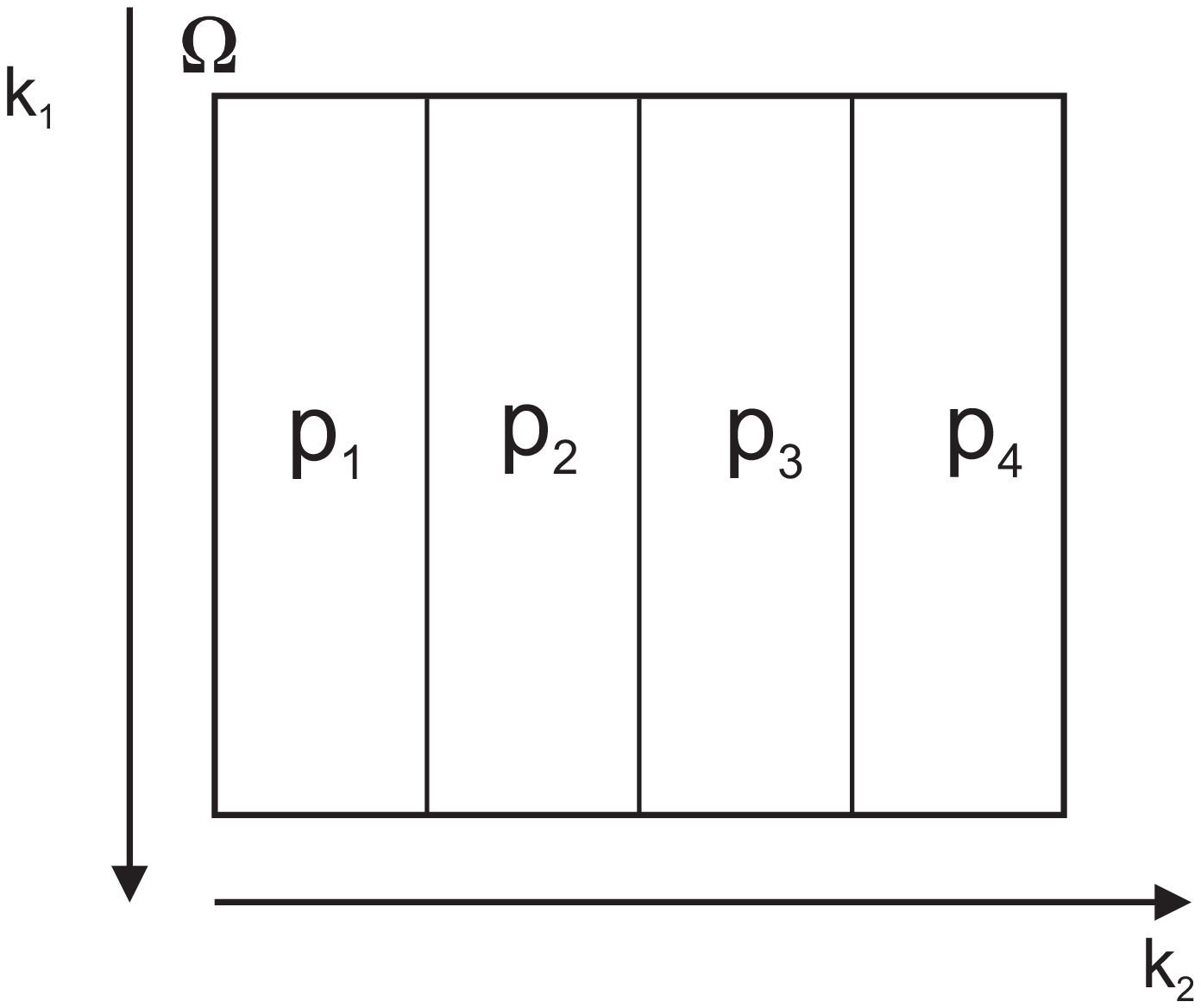}\hfill
\includegraphics[width=0.5\textwidth]{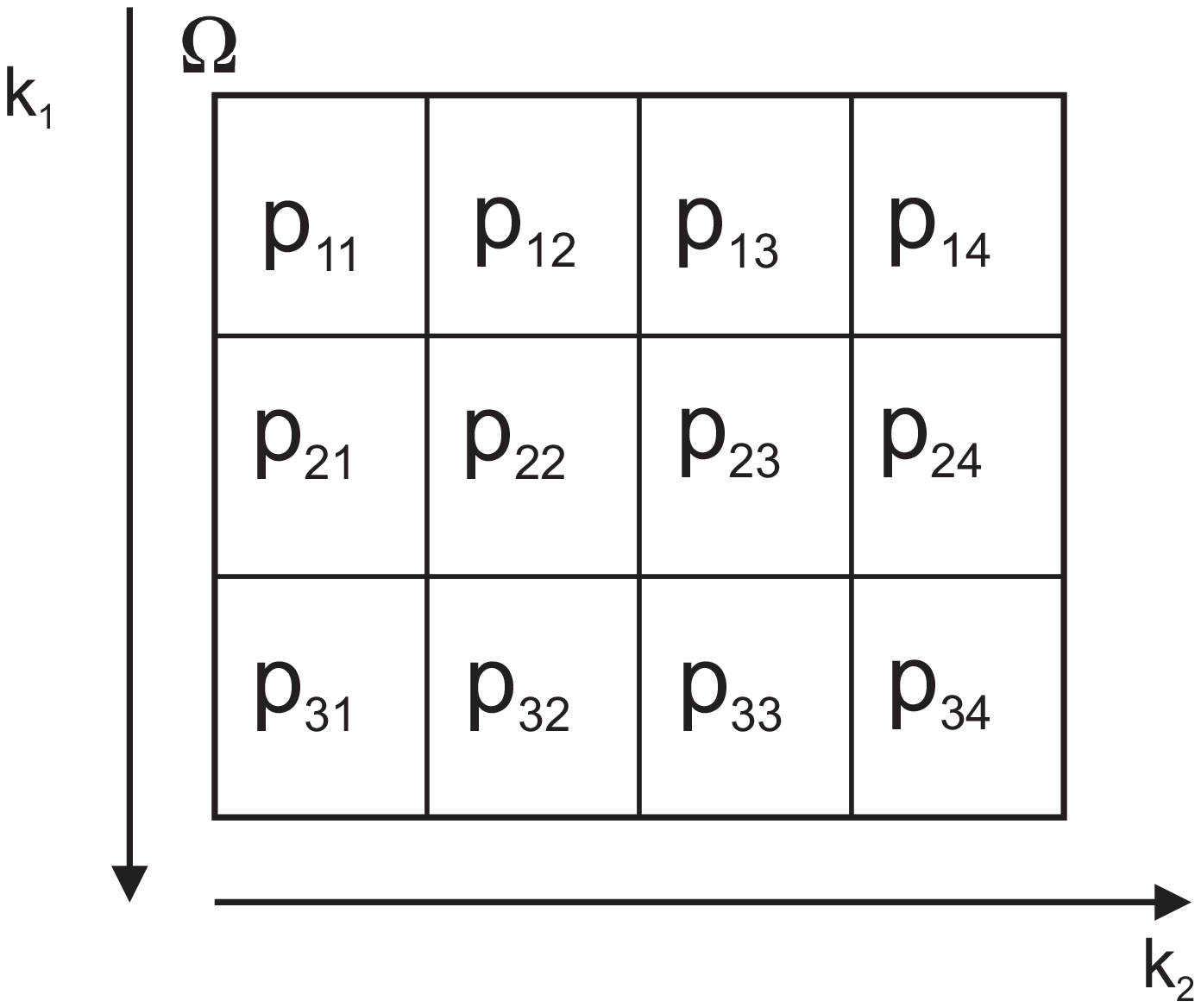}\\
\parbox{0.5\textwidth}{\center a) Method of variable separation}\hfill
\parbox{0.5\textwidth}{\center b) Method of alternating directions }
\caption{Domain Decomposition}
\end{center}
\end{figure}

Test calculations were performed on an MBC-100k supercomputer of
the Interdepartment Center of the Russian Academy of Sciences; the
supercomputer is based on Intel Xeon four-core processors
operating at 3 GHz in the Infiniband communication environment.

Results obtained for the dependence of the computing time
(\textbf{$T_{avr}$}) and speedup (\textbf{$S_{avr}$}) for the
Fourier and ADI methods are listed in Tables $1$ and $2$, and in
(Figs. 5a,5b,6a,6b).

Based on the obtained results, we will point out the following:

\begin{itemize}
    \item For the Fourier and ADI methods, the dependence of the computing time on the number of processors is linear.
    \item For computing by the ADI method, starting from some number of processors, the speedup is superlinear because
    as the number of processors grows, the data volume per PE decreases, therefore,
    they can be located completely in a faster memory cache.
    \item The maximum performance of the Fourier method was $1700$ equations/ sec. for a $512\mathrm{x}512$ mesh $833$ eqs./sec for $1024\mathrm{x}1024$, $417$ eqs./sec for $2048 2048$,
     $161$  eqs./sec for $4096\mathrm{x} 4096$, $56$ eqs./sec for $8192\mathrm{x}8192$, and $13$ eqs. for $16384\mathrm{x}16384$, respectively.
    \item When the number of nodes in
    one direction exceeds several times the number of nodes
    in another direction, the efficiency of the parallel Fourier algorithm is between $80\%$ and $95\%$ (Figs. 6a and 6b).
\end{itemize}

\begin{table}[!h]
\center \small
\begin{tabular}{|l|c|c||c|c||c|c||c|c||c|c||c|c|}

  \hline
  size & \multicolumn{2}{c||}{512x512}& \multicolumn{2}{c||}{1024x1024}&\multicolumn{2}{c||}{2048x2048}&\multicolumn{2}{c||}{4096x4096}&\multicolumn{2}{c||}{8192x8192}&\multicolumn{2}{c|}{16384x16384}   \\ \hline
   NP &  $ \mathrm{T}_{avr}$&  $ \mathrm{S}_{avr}$  & $ \mathrm{T}_{avr}$&$ \mathrm{S}_{avr}$&$ \mathrm{T}_{avr}$&$ \mathrm{S}_{avr}$&$ \mathrm{T}_{avr}$&$ \mathrm{S}_{avr}$&$ \mathrm{T}_{avr}$&$ \mathrm{S}_{avr}$&$ \mathrm{T}_{avr}$&$ \mathrm{S}_{avr}$\\
  \hline

  1 &  2.6e-02 & -&1.1e-01& -&  4.9e-01&-&2.11&-&   11 &- &- & -\\ \hline
  4 &   6.2e-03& 4.2&2.5e-02&4.5&1.2e-01&4& 5.5e-01& 3.8& 2.8  & 3.9 &-&- \\ \hline
  8 &  2.6e-03& 10&1.1e-02&10&6e-02&8.2&2.9e-01&7.2&  1.6  &6.7& - &- \\\hline
  16 & 1.4e-03 & 18.5& 5.6e-03&21&3e-02&16.3&1.3e-01&   16.2&  0.8 & 13.8&-& - \\ \hline
  32 &  8.5e-04 &30&3e-03 & 38&1.3e-02&38& 6.7e-02& 31.4 & 0.4 & 27.5&1.78 &-\\ \hline
  64 &  7.9e-04&  33&2e-03 &58&6.6e-03&74&3.3e-02&63.9& 0.19 & 58.3&0.95&54 \\ \hline
  128 &   5.9e-04& 44&1.3e-03&84& 4e-03&122&  1.5e-02&140&  9.6e-02&115.4& 0.45& 126 \\ \hline
  256 &  1.2e-03 &22&1.2e-03 & 96&2.8e-03&175&9.4e-03&224&  5e-02&221.6&0.24&237 \\ \hline
  512 & - &-& 2.2e-03 & 52&2.4e-03&204&  6.8e-03&310&  2.8e-02& 395&1.4e-01&406 \\ \hline
  1024 & - &-&-&-&-&-&   6.2e-03&  340&  1.8e-02&611&7.7e-02&   739\\\hline

  \hline

\end{tabular}
\caption{ Computing time ($\mathrm{T}_{avr}$) and speedup
($S_{avr}$) versus the number of processors for the Fourier
method.}
\end{table}

\begin{table}[!h]
\center
\small
\begin{tabular}{|l|c|c|c||c|c|c||c|c|c||c|c|c|}

  \hline
  size & \multicolumn{3}{c||}{512x512}& \multicolumn{3}{c||}{1024x1024}&\multicolumn{3}{c||}{2048x2048}&\multicolumn{3}{c|}{4096x4096} \\ \hline

  NP & $ \mathrm{T}_{avr}$&  $\mathrm{S}_{avr}$ &  $\mathrm{M}$ &$ \mathrm{T}_{avr}$&$\mathrm{S}_{avr}$&$\mathrm{M}$&$\mathrm{T}_{avr}$&$\mathrm{S}_{avr}$&$\mathrm{M}$&$\mathrm{T}_{avr}$&$\mathrm{S}_{avr}$&$\mathrm{M}$\\
  \hline

  1 &  0.9 &  -&- &8.7&-&-&48.17 &  -&- &202&-&- \\ \hline
  4 &  8.5e-02&  10.5& 2x2&1.1 &7.9&1x4&9.6&  5 & 1x4&67&  3& 1x4\\ \hline
  8 &  6.7e-02& 13.4&  8x1 &0.84&10.3&1x8&7.3& 6.6&  1x8&34& 5.9&  1x8\\\hline
  16 & 2.6e-02 & 34 &  1x16&2.8e-01&31&1x16& 2.3 & 21 &  1x16&12.2 & 16.8 &1x16\\ \hline
  32 &  2.0e-02 &45  &  2x16 &   8.1e-02  &107&2x16&  0.94  & 51   &  2x16&7.7   & 28 &  1x32\\ \hline
  64 &  2.5e-02&  36&  4x16 & 4.6e-02&189&4x16 &  0.29& 166   &  4x16 &2.8&  71&  1x64\\ \hline
  128 &  -& -  & -&3.7e-02& 235&  4x32&  8.3e-02&  580&  4x32&1&  180&  4x32\\ \hline
  256 & - &- &-&3.4e-02 & 255&  16x16 & 5e-02 &  963&  16x16&0.3&  721&  16x16\\ \hline
  512 & - &-& -&3.3e-02 &263& 16x32&  4.6e-02& 1047&  16x32&9.7e-02&  2082&  16x32  \\ \hline
  1024 & - &-   &-& 4.2e-02 &  207 &  32x32&  4.7e-02& 1024  &  32x32&6.8e-02&  2970 & 32x32\\

  \hline

\end{tabular}
\caption{Computing time ($\mathrm{T}_{avr})$ and speedup
($S_{avr}$) versus the number of processors for the ADI method.
($M$, the number of processors in directions $k_1$ and $k_2$,
which enabled the minimal computing time.)}
\end{table}

\begin{figure}[!h]
\label{main_pic1}
\begin{center}
\includegraphics[width=0.5\textwidth]{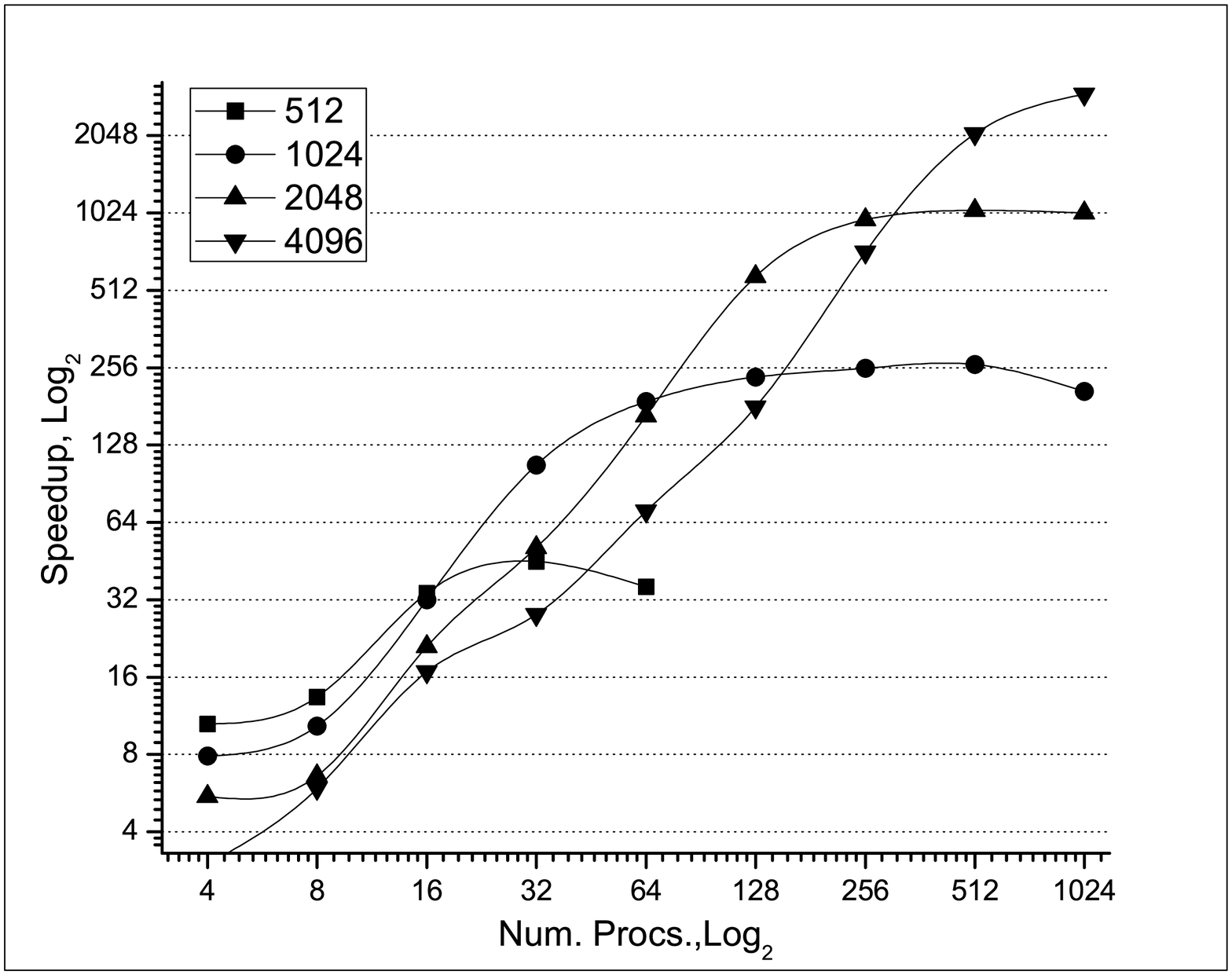}\hfill
\includegraphics[width=0.5\textwidth]{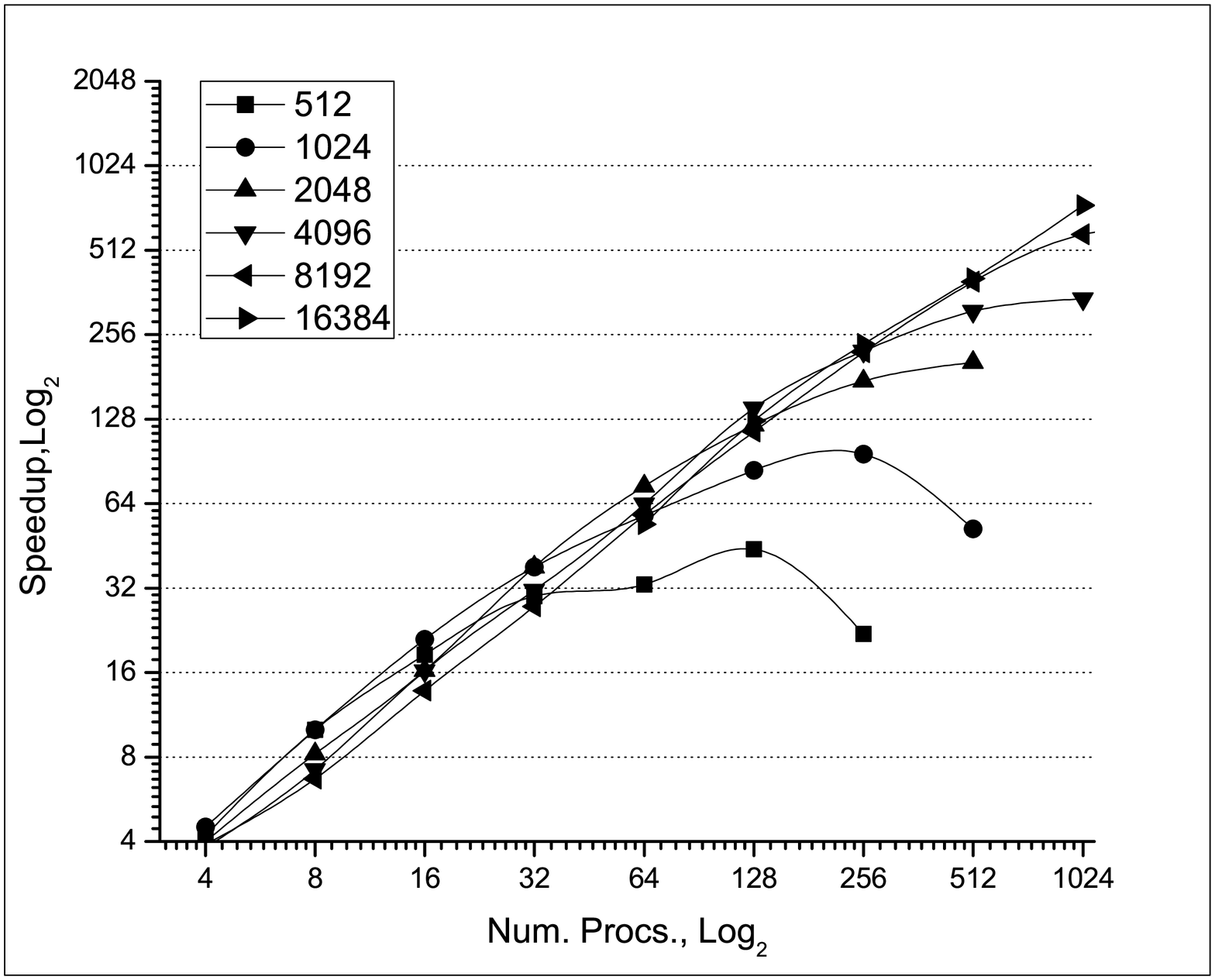} \\
\parbox{0.5\textwidth}{\center a) Method of alternating directions}\hfill
\parbox{0.5\textwidth}{\center b) Method of variable separation }
\caption{Speedup versus the number of processors for different
meshes}
\end{center}
\end{figure}

\begin{figure}[!h]
\label{main_pic2}
\includegraphics[width=0.5\textwidth]{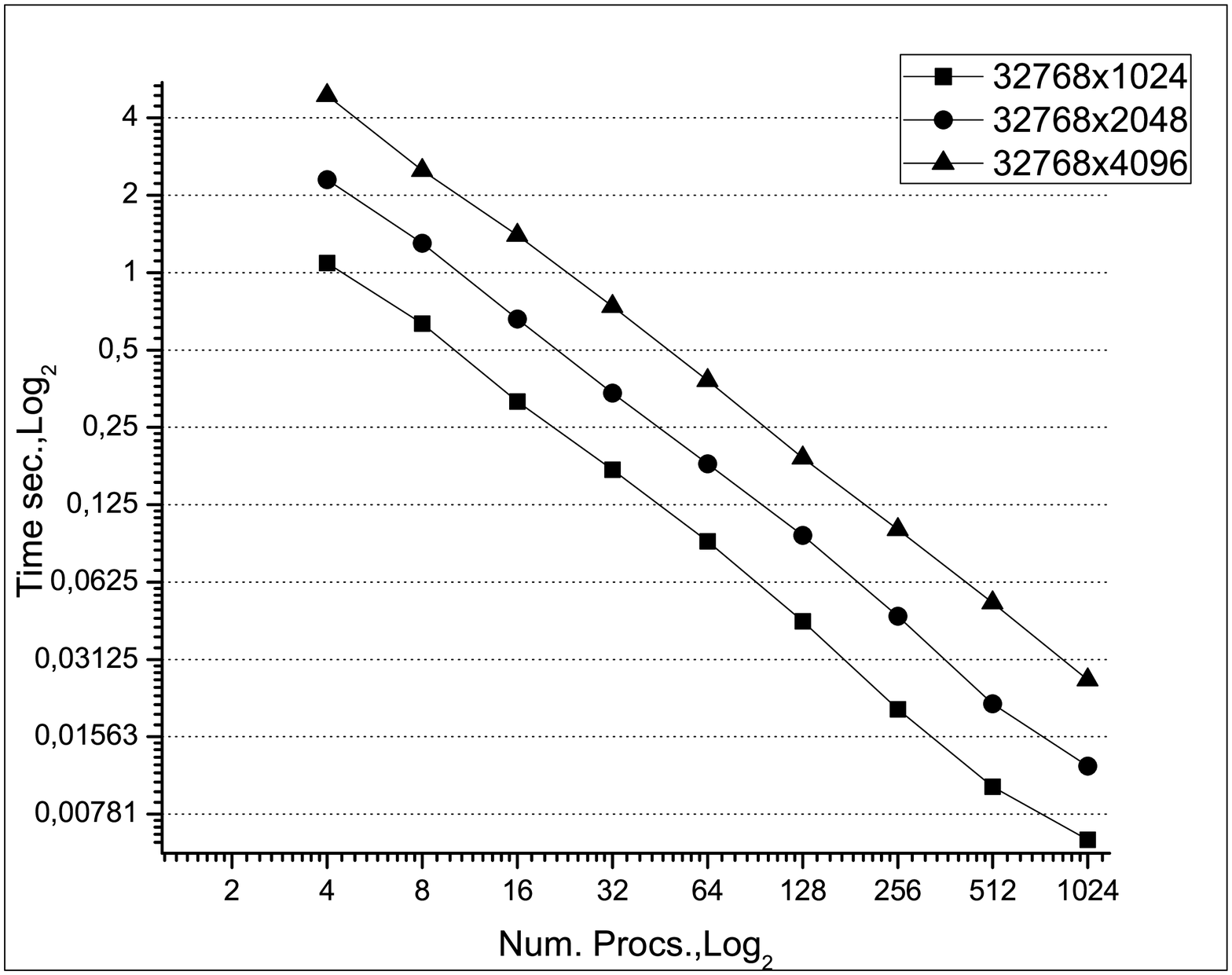} \hfill
\includegraphics[width=0.5\textwidth]{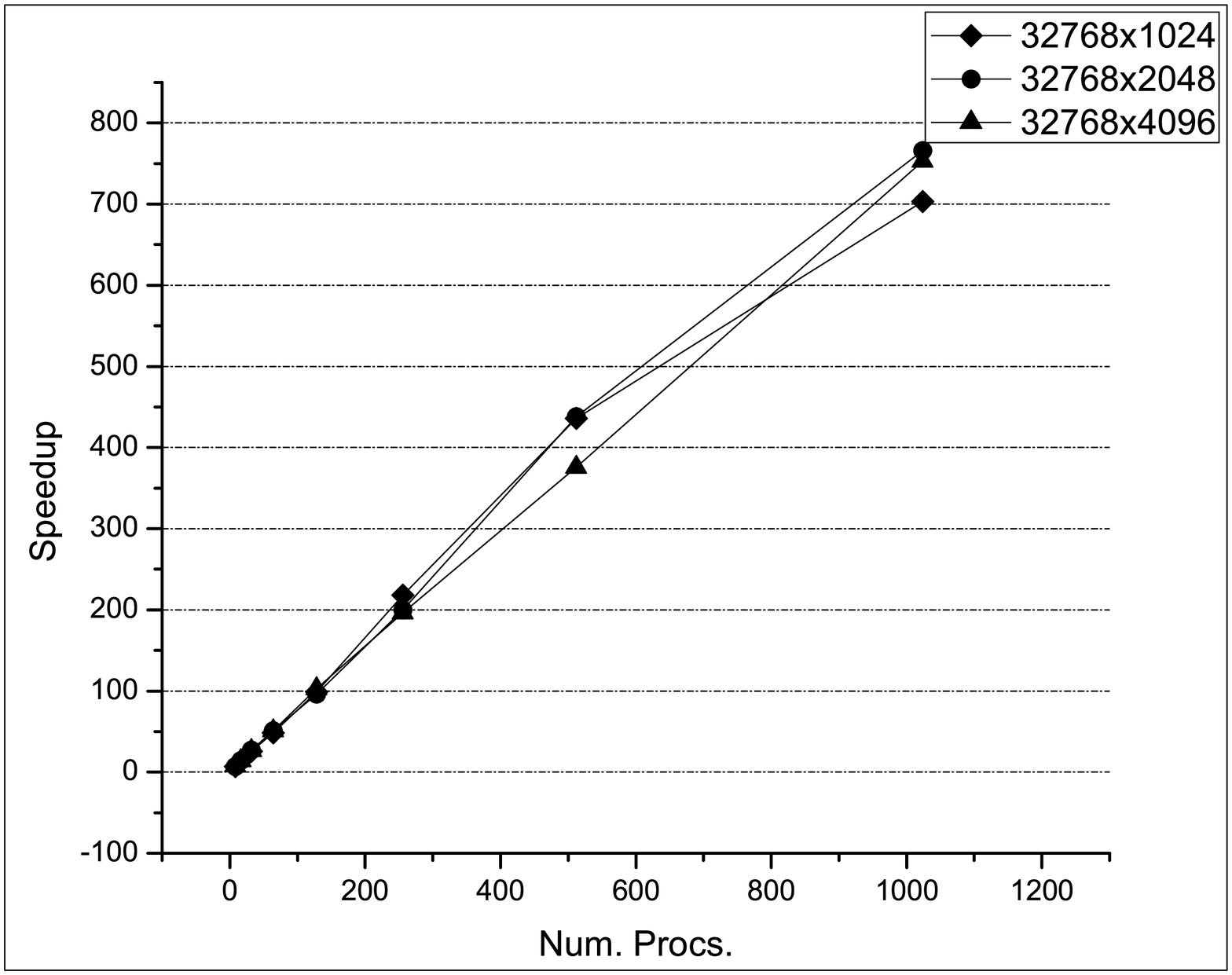}\\
\parbox{0.5\textwidth}{\center a)}\hfill
\parbox{0.5\textwidth}{\center b)}

\caption{Computing time (a) and speedup (b) for the method of
variable separation in the case of a rectangular region versus the
number of processors}
\end{figure}

Thus, as a result of our computational experiments we registered
the presently record efficiency in solving the Poisson equation on
a multicomputer. These results were achieved due to applying the
dichotomy algorithm for a series of tridiagonal SLAEs, which was
specially designed for distributed-memory supercomputers. We
should note that the proposed algorithms will be no less
efficient, but even more for implementing on shared-memory
multiprocessor computer systems because the communication
interactions are minimal in that case.

\section{Conclusions}
The proposed parallel sweep algorithm for solving a series of
tridiagonal systems of linear algebraic equations has validated
its efficiency as a result of computational experiments. The main
feature of the algorithm is that it is required at first to
perform some preliminary computations whose complexity is
comparable with solving one problem, and then solve a number of
SLAEs for different right-hand sides with a nearly linear speedup.
Thus, we have developed and investigated a promising method for
solving a series of tridiagonal SLAEs, whose efficiency and
scalability are record for today.

\section{Acknowledgments.} The author is grateful to Prof. V. Malyshkin
for fruitful discussions of the paper.

\end{document}